\documentclass[11pt]{article}
\usepackage{amsfonts,latexsym,rawfonts,amsmath,amssymb,amsthm,graphicx}

\numberwithin{equation}{section}
\newtheorem{theorem}{Theorem}[section]
\newtheorem{lem}[theorem]{Lemma}
\newtheorem{thm}[theorem]{Theorem}
\newtheorem{pro}[theorem]{Proposition}
\newtheorem{cor}[theorem]{Corollary}
\newtheorem{defi}[theorem]{Definition}
\newtheorem{rem}[theorem]{Remark}

\def\dist{{\rm dist\,}}

\def\endproof{$\hfill\Box$\\}
\def\N{\mathbb{N}}
\def\R{\mathbb{R}}
\def\Z{\mathbb{Z}}

\def\C{\mathbb{C}}
\def\W{\mathcal{W}}
\def\CI{W^{2,2}_{{\rm conf}}}
\def\genus{p}
\def\diam{{\rm diam\,}}
\def\arccot {{\rm arccot}}

\title{\bf $W^{2,2}$-conformal immersions of
a closed Riemann surface into $\R^n$}
\author{{\sc Ernst Kuwert} and {\sc Yuxiang Li}
\footnote{The authors were supported by DFG Collaborative Research Center SFB/Transregio 71.}
\footnote{The second author was partially supported by the National Science Foundation of China (No. 10801082)
and the Specialized Research Fund for the Doctoral Program of Higher Education (No. 200800031078).}}
\date{}
\begin{document}
\maketitle

\begin{abstract}
{\noindent
We study sequences $f_k:\Sigma_k \to \R^n$ of conformally
immersed, compact Riemann surfaces with fixed genus and
Willmore energy ${\cal W}(f) \leq \Lambda$. Assume that
$\Sigma_k$ converges to $\Sigma$ in moduli space, i.e.
$\phi_k^\ast(\Sigma_k) \to \Sigma$ as complex structures
for diffeomorphisms $\phi_k$. Then we construct
a branched conformal immersion $f:\Sigma \to \R^n$ and
M\"obius transformations $\sigma_k$, such that for a
subsequence $\sigma_k \circ f_k \circ \phi_k \to f$
weakly in $W^{2,2}_{loc}$ away from finitely many points.
For $\Lambda < 8\pi$ the map $f$ is unbranched.
If the $\Sigma_k$ diverge in moduli space, then we show
$\liminf_{k \to \infty} {\cal W}(f_k) \geq \min(8\pi,\omega^n_p)$.
Our work generalizes results in \cite{K-S3} to arbitrary
codimension.}
\end{abstract}

\section{Introduction}
Let $\Sigma$ be a closed oriented surface of genus $\genus \in \N_0$.
For an immersion $f:\Sigma \to \R^n$ the Willmore functional is
defined by
$$
\W(f)=\frac{1}{4}\int_\Sigma|H|^2\,d\mu_g,
$$
where $H$ is the mean curvature vector and $g$ is the induced
metric on $\Sigma$. The infimum among closed immersed surfaces
of genus $p$ is denoted by $\beta^n_p$.
We have $\beta^n_0 = 4\pi$, which is attained only by round
spheres \cite{W}. For $p \geq 1$ we have the inequalities
$4\pi < \beta^n_p < 8\pi$ \cite{S,K}.
In this paper we study compactness properties of sequences
$f_k:\Sigma \to \R^n$ with ${\cal W}(f_k) \leq \Lambda$.
By the Gau{\ss} equations and Gau{\ss}-Bonnet, the second
fundamental form is then equivalently bounded by
$$
\int_{\Sigma} |A_{f_k}|^2\,d\mu_{g_k} \leq 4 \Lambda + 8\pi(p-1).
$$
In \cite{L} Langer proved a compactness theorem for surfaces with
$\|A\|_{L^q} \leq \Lambda$ for $q > 2$, using that the surfaces
are represented as $C^1$-bounded graphs over disks of radius
$r(n,q,\Lambda) > 0$. Clearly, the relevant Sobolev embedding
fails for $q = 2$. For surfaces with $\|A\|_{L^2}$ small in a
ball, L. Simon proved an approximate graphical decomposition,
see \cite{S}, and showed the existence of Willmore minimizers
for any $p \geq 1$, assuming for $p \geq 2$ that
$$
\beta^n_p <
\min\Big\{4\pi+\sum\limits_{i}(\beta_{\genus_i}^n-4\pi):
\sum\limits_{i} \genus_i = \genus,\,1 \leq p_i < p\Big\} = \omega^n_p.
$$
This inequality was confirmed later in \cite{B-K}. As
$\lim_{p \to \infty} \beta^n_p = 8\pi$ by  \cite{K-L-S},
we have $\omega^n_p > 8\pi$ for large $p$.
Recently, using the annulus version of the approximate
graphical decomposition lemma, a compactness theorem was proved
in \cite{K-S3} for surfaces in $\R^3$ under the assumptions
$$
\liminf_{k \to \infty}  {\cal W}(f_k) <
\begin{cases}
8 \pi & \mbox{ if } p = 1,\\
\min(8\pi,\omega^3_p) & \mbox{ if } p \geq 2.
\end{cases}
$$
Moreover, it was shown that these conditions are optimal. For
$n = 4$ the result was proved under the additional assumption
$\liminf_{k \to \infty} {\cal W}(f_k) < \beta^4_p + \frac{8\pi}{3}$.
In \cite{K-S4} these compactness theorems were applied to prove the
existence of a Willmore minimizer with prescribed conformal type.\\
\\
Here we develop a new approach to compactness, generalizing
the results of \cite{K-S3} to any codimension. As main
tools we use a convergence theorem of Hurwitz type for conformal
immersions, which is due to H\'elein \cite{H}, and the estimates
for the conformal factor by M\"uller and \v{S}ver\'ak \cite{M-S}.
The paper is organized as follows. In Section 2 we introduce
the notion of $W^{2,2}$ conformal immersions, and recall the
main estimate from \cite{M-S} as well as the monotonicity formula
from \cite{S}. In Section 3 we adapt the analysis of \cite{M-S}
to show that isolated singularities of conformal immersions
with square integrable second fundamental form and finite area
are branchpoints, in a suitable weak sense. The compactness theorem
for conformal immersions is presented in Section 4. We first
deal with the case of a fixed Riemann surface in Proposition
\ref{fixcc}, and extend the result to sequences of Riemann surfaces
converging in moduli space in Theorem \ref{convergence2}.
Finally in Section 5, we study surfaces whose conformal type
degenerates and show that the lower bound from \cite{K-S3}
extends to higher codimension. Along the lines, we state a
version of Theorem 5.1.1 in \cite{H} with optimal constants.

\section{$W^{2,2}$ conformal immersions}

\begin{defi} \label{defconformalimmersion}
Let $\Sigma$ be a Riemann surface. A map $f\in W^{2,2}_{{\rm loc}}(\Sigma,\mathbb{R}^n)$
is called a conformal immersion, if in any local parameter the induced metric
$g_{ij} = \langle \partial_i f,\partial_j f \rangle$ is given by
$$
g_{ij} = e^{2u} \delta_{ij} \quad \mbox{ where } u \in L^\infty_{{\rm loc}}(U).
$$
For compact $\Sigma$ the set of all $W^{2,2}$-conformal immersions  is denoted by
$\CI(\Sigma,\R^n)$.
\end{defi}
It is easy to see that for $f \in \CI(\Sigma,\R^n)$ one has in a local parameter
$$
\textstyle{u = \frac{1}{2} \log \left(\frac{1}{2} |Df|^2\right) \in W^{1,2}_{{\rm loc}}(U).}
$$
The induced measure $\mu_g$, the second fundamental form $A$ and the
mean curvature vector $H$ are given by the standard coordinate
formulas. We define $K_g$ by the Gau{\ss} equation
$$
K_g = \frac{1}{2} (|H|^2 - |A|_g^2)
= e^{-4u} (\langle A_{11},A_{22}\rangle - |A_{12}|^2).
$$
In a local parameter, we will now verify the weak Liouville equation
$$
\int_U \langle Du, D\varphi \rangle = \int_U K_g e^{2u}\varphi \quad
\mbox{ for all }\varphi \in C^\infty_0(U).
$$
In particular, this shows that $K_g$ is intrinsic. We start by computing
$$
\langle \partial^2_{ij} f, \partial_k f \rangle
+ \langle \partial^2_{ki} f, \partial_j f \rangle
= 2 e^{2u} \partial_i u\, \delta_{jk},
$$
which implies after permutation of the indices that
$$
\langle \partial^2_{ij} f,\partial_k f \rangle =
e^{2u}\big(\partial_i u \, \delta_{jk} + \partial_j u \,\delta_{ik} - \partial_k u\, \delta_{ij}\big).
$$
Expanding explicitely yields
\begin{eqnarray*}
\partial^2_{11} f & = & A_{11} + \partial_1 u\, \partial_1 f - \partial_2 u\, \partial_2 f\\
\partial^2_{22} f & = & A_{22} - \partial_1 u\, \partial_1 f + \partial_2 u\, \partial_2 f\\
\partial^2_{12} f & = & A_{12} + \partial_2 u\, \partial_1 f + \partial_1 u\, \partial_2 f,
\end{eqnarray*}
and we get
$$
\langle A_{11},A_{22}\rangle - |A_{12}|^2 =
\langle \partial^2_{11} f, \partial^2_{22} f \rangle - |\partial^2_{12} f|^2
+ 2 e^{2u} |Du|^2.
$$
For any $u \in  W^{1,2} \cap L^{\infty}(U)$, $f \in W^{2,2}_{{\rm loc}}(U,\R^n)$
and $\varphi \in C^\infty_0(U)$ we have the formula
$$
\int_U \Big(\langle \partial^2_{11} f, \partial^2_{22} f \rangle - |\partial^2_{12} f|^2\Big)
e^{-2u} \varphi =
\int_U \Big(\langle \partial_1 f,\partial^2_{12} f \rangle \partial_2 \big(e^{-2u}\varphi\big)
- \langle \partial_1 f,\partial^2_{22} f \rangle \partial_1 \big(e^{-2u}\varphi\big)\Big).
$$
This follows by approximation from the case when $f$ is smooth.
Now for $f$ conformal
$$
\langle \partial_1 f,\partial^2_{12} f \rangle = e^{2u} \partial_2 u
\quad \mbox{ and } \quad
\langle \partial_1 f,\partial^2_{22} f \rangle  =  - e^{2u} \partial_1 u,
$$
which yields
$$
\int_U \Big(\langle \partial^2_{11} f, \partial^2_{22} f \rangle - |\partial^2_{12} f|^2\Big)
e^{-2u} \varphi = \int_U \langle Du,D\varphi \rangle - 2 \int_U |Du|^2 \varphi,
$$
and the Liouville equation follows.

\begin{rem}
More generally if $g = e^{2u} g_0$ where $g_0$ is any smooth conformal
metric, then
$$
-\Delta_{g_0} u = K_g e^{2u} - K_{g_0} \quad \mbox{ weakly.}
$$
Testing with a constant function, we infer for
closed $\Sigma$ the Gau{\ss}-Bonnet formula
$$
\int_\Sigma K_g\,d\mu_g = 2\pi\chi(\Sigma).
$$
\end{rem}

$W^{2,2}$ conformal immersions $f$ can be
approximated by smooth immersions in the $W^{2,2}$ norm.
In fact, a standard mollification $f_\varepsilon$ will be
immersed for small $\varepsilon > 0$, by an argument of
\cite{S-U}.

\subsection{Gau{\ss} map and compensated compactness}

By assumption the right hand side $K_g e^{2u}$ of the Liouville
equation belongs to $L^1$. In \cite{M-S} M\"uller and
\v{S}ver\'ak discovered that the term can be written
as a sum of Jacobi determinants, and that improved estimates
can be obtained from the Wente lemma \cite{Wen} or from
\cite{C-L-M-S}. The following result is Corollary 3.5.7
of \cite{M-S}. Recall that in their notation,
$\omega$ denotes twice the standard K\"ahler
form and $W^{1,2}_0(\C)$ is the space of functions
$v \in L^2_{loc}(\C)$ with $Dv\in L^2(\C)$.

\begin{thm}\label{MS}
Let $\varphi \in W^{1,2}_0(\C,\C P^n)$ satisfy
$$
\int_{\C} \varphi^{*}\omega =0 \quad \mbox{ and } \quad
\int_{\C} J\varphi \leq \gamma < 2\pi.
$$
Then there is a unique function $v \in W^{1,2}_0(\C)$
solving the equation $- \Delta v = \ast \varphi^{\ast} \omega$ in $\C$
with boundary condition $\lim_{z \to \infty} v(z) = 0$. Moreover
$$
\|v\|_{L^\infty(\C)} +\|D v\|_{L^2(\C)}
\leq C(\gamma) \int_{\C} |D\varphi|^2.
$$
\end{thm}

For $f\in W^{2,2}_{{\rm conf}}(D,\R^n)$ let $G \in W^{1,2}(D,\C P^{n-1})$
be the associated Gau{\ss} map. Here we embed the Grassmannian
$G(2,n)$ of oriented $2$-planes into $\C P^{n-1}$ by sending an
orthonormal basis $e_{1,2}$ to $[(e_1+i e_2)/\sqrt{2}]$. Then
$$
K_g e^{2u}  = \ast G^{\ast} \omega \quad \mbox{ and } \quad
\int_D |DG|^2 = \frac{1}{2} \int_D |A|^2\,d\mu_g.
$$
\begin{cor}\label{key} For $f\in W^{2,2}_{{\rm conf}}(D,\R^n)$
with induced metric $g_{ij} = e^{2u} \delta_{ij}$, assume
$$
\int_D |A|^2\,d\mu_g \leq \gamma < \gamma_n =
\begin{cases}
8\pi & \mbox{ if } n= 3,\\4\pi & \mbox{ if }n \geq 4.
\end{cases}
$$
Then there exists a  function $v:\C \to \R$ solving the equation
$$
-\Delta v = K_g e^{2u} \quad \mbox{ in } D,
$$
and satisfying the estimates
$$
\|v\|_{L^\infty(\C)} + \|Dv\|_{L^2(\C)}
\leq C(\gamma) \int_D |A|^2\,d\mu_g.
$$
\end{cor}

\proof

We follow \cite{M-S}. Define the map $\varphi:\C \to \C P^{n-1}$ by
$$
\varphi(z) =
\left\{\begin{array}{ll}
G(z) & \mbox{ if }z \in D\\
G(\frac{1}{\overline{z}}) & \mbox{ if }z \in \C \backslash \overline{D}.
\end{array}\right.$$
Then $\varphi\in W_0^{1,2}(\mathbb{C},\C P^{n-1})$ and
$\int_{\mathbb{C}}{\varphi}^* \omega =0$. For $n \geq 4$ we have
$$
\int_{\mathbb{C}} J\varphi = 2 \int_D JG \leq
\frac{1}{2} \int_D |A|^2\,d\mu_g \leq \frac{\gamma}{2} < 2\pi.
$$
Thus the result follows from Theorem \ref{MS}. The same is
true for $n=3$, since then
$$
\int_{\C} J\varphi = \frac{1}{2} \int_D |K_g|\,d\mu_g \leq
\frac{1}{4} \int_D |A|^2\,d\mu_g \leq \frac{\gamma}{4} < 2\pi.
$$
\endproof

The function $K_g e^{2u}$ belongs actually to the Hardy space
$\mathcal{H}^1$, see \cite{C-L-M-S}. This implies that $v$ has
second derivatives in $L^1$, and in particular that $v$ is
continuous \cite{M-S}. As $u-v$ is harmonic, it follows that $u$
is also continuous, but this will not be used here. The following
is an immediate consequence of Corollary \ref{key}.

\begin{cor}\label{key2} Let $f \in W^{2,2}_{{\rm conf}}(D,\R^n)$
with induced metric $g_{ij} = e^{2u} \delta_{ij}$. If
$$
\int_D |A|^2\,d\mu_g \leq \gamma < \gamma_n,
$$
then we have the estimate
$$
\|u\|_{L^\infty(D_\frac{1}{2})} + \|Du\|_{L^2(D_\frac{1}{2})}
\leq C(\gamma) \Big(\int_D |A|^2\,d\mu_g + \|u\|_{L^1(D)}\Big).
$$
\end{cor}


\subsection{Simon's monotonicity formula}
We briefly review the monotonicity identity from \cite{S}
for proper $W^{2,2}$ conformal immersions $f:\Sigma \to \R^n$.
For more details we refer to \cite{K-S2}. Since $f$ is locally
Lipschitz, the measure $\mu = f(\mu_g)$ is an integral varifold
with multiplicity function $\theta^2(\mu,x) = \# f^{-1}\{x\}$
and approximate tangent space $T_x \mu = Df(p)\cdot T_p\Sigma$
a.e. when $x = f(p)$. The immersion $f$ satisfies
$$
\label{eqci2}
\int_\Sigma {\rm div}_g X \,d\mu_g =
- \int_\Sigma \langle X, H \rangle\,d\mu_g \quad
\mbox{ for any }X \in W^{1,1}_0(\Sigma,\R^n).
$$
For the varifold $\mu$ this implies the first variation formula
$$
\int_U {\rm div}_{\mu} \phi\,d\mu = - \int_U \langle \phi,H_\mu \rangle\,d\mu,
$$
where the weak mean curvature is given by
$$
H_\mu(x) = \begin{cases}
\frac{1}{\theta^2(\mu,x)} \sum_{p \in f^{-1}\{x\}} H(p) &
\mbox{ if } \theta^2(\mu,x) > 0,\\
0 & \mbox{ else. }
\end{cases}
$$
From the definiton we have trivially the inequality
$$
{\cal W}(\mu,V) =
\frac{1}{4} \int_V |H_\mu|^2\,d\mu
\leq \frac{1}{4} \int_{f^{-1}(V)} |H|^2\,d\mu_g.
$$
Observing that $H_\mu(x)$ is $\mu$ a.e. perpendicular to $T_x\mu$, the
proof of the monotonicity identity in \cite{S} extends to show
that for $B_\sigma(x_0) \subset B_\varrho(x_0)$ one has
$$
g_{x_0}(\varrho) - g_{x_0}(\sigma) =
\frac{1}{16\pi} \int_{B_\varrho(x_0) \backslash B_\sigma(x_0)}
\Big|H_\mu + 4\frac{(x-x_0)^\perp}{|x-x_0|^2}\Big|^2\,d\mu,
\quad \mbox{ where }
$$
$$
g_{x_0}(r) =
\frac{\mu(B_r(x_0))}{\pi r^2}
+ \frac{1}{4 \pi} {\cal W}(\mu,B_r(x_0))
+ \frac{1}{2\pi r^2} \int_{B_r(x_0)} \langle x-x_0,H_\mu \rangle \,d\mu.
$$
Applications include the existence and upper semicontinuity of $\theta^2(\mu,x)$
and, for closed surfaces, the Li-Yau inequality, see \cite{LY},
$$
\theta^2(\mu,x) \leq \frac{1}{4\pi} {\cal W}(f).
$$
Another consequence is the diameter bound from \cite{S}. If $\Sigma$ is compact
and connected, then for $f \in \CI(\Sigma,\R^n)$ one obtains
\begin{equation}\label{diameter}
\Big(\frac{\mu_g(\Sigma)}{{\cal W}(f)}\Big)^{\frac{1}{2}}
\leq \diam f(\Sigma) \leq
C \Big(\mu_g(\Sigma)\,{\cal W}(f)\Big)^{\frac{1}{2}}.
\end{equation}

\section{Classification of isolated singularities}
In \cite{M-S} M\"uller and \v{S}ver\'{a}k studied
the behavior at infinity of complete, conformally
parametrized surfaces with square integrable second
fundamental form.  Here we adapt their analysis to
the case of finite isolated singularities.

\begin{thm}\label{removal}
Suppose that $f\in W^{2,2}_{{\rm conf},loc}(D\backslash \{0\},\R^n)$ satisfies
$$
\int_D |A|^2\,d\mu_g < \infty \quad \mbox{ and } \quad \mu_g(D) < \infty,
$$
where $g_{ij} = e^{2u} \delta_{ij}$ is the induced metric. Then
$f \in W^{2,2}(D,\R^n)$ and we have
\begin{eqnarray*}
u(z) & = & m\log |z|+ \omega(z) \quad \mbox{ where }
m \in \N_0,\,\omega \in C^0 \cap W^{1,2}(D),\\
-\Delta u & = & -2m\pi \delta_0+K_g e^{2u} \quad \mbox{ in }D.
\end{eqnarray*}
The multiplicity of the immersion at $f(0)$ is given by
$$
\theta^2\big(f(\mu_g \llcorner D_\sigma(0)),f(0)\big) = m+1 \quad \mbox{ for any small }
\sigma > 0.
$$
Moreover, if $m = 0$ then $f$ is a conformal immersion on $D$.
\end{thm}

\proof

We may assume $\int_D |A|^2 \,d\mu_g < 4\pi$, hence the
Gau{\ss} map $G:D \to G(n,2)$ has energy
$$
\int_D |DG|^2 = \frac{1}{2} \int_D |A|^2\,d\mu_g < 2\pi.
$$
Extending by $G(z):=G(1/\overline{z})$ for
$|z| > 1$ yields $G \in W^{1,2}_0\big(\R^2,G(n,2)\big)$,
where
$$
\int_{\R^2} G^\ast \omega = 0 \quad \mbox{ and } \quad
\int_{\R^2} JG \leq \int_D |DG|^2 < 2\pi.
$$
Thus there is a function $v  \in C^0 \cap W^{1,2}_0(\R^2)$
such that
$$
- \Delta v = K_g e^{2u} \quad \mbox{ and } \quad\lim_{z \to \infty} v(z) = 0,
$$
$$
\|v\|_{C^0(\R^2)} + \|Dv\|_{L^2(\R^2)} \leq C \int_D |A|^2\,d\mu_g.
$$
Now consider the harmonic function
$h:D\backslash \{0\} \to \R$,\,$h(z) = u(z)-v(z) - \alpha \log |z|$, where
$$
\alpha = \frac{1}{2\pi} \int_{\partial D_r(0)}
\frac{\partial (u-v)}{\partial r}\,ds \in \R
\quad \mbox{ for }r \in (0,1).
$$
We claim that $h$ has a removable singularity at the origin.
Let $h = {\rm Re\,}\phi$ where $\phi:D \backslash \{0\} \to \C$ is
holomorphic, and compute for $m = \min \{k \in \Z: k \geq \alpha\}$
$$
\big|z^m e^{\phi(z)}\big| = |z|^m e^{h(z)} \leq
e^{u(z)-v(z)} \leq C e^{u(z)} \in L^2(D).
$$
Thus $z^m e^{\phi(z)} = z^k g(z)$ for $k \in \N_0$ and
$g:D \to \C\backslash \{0\}$ holomorphic, which yields
$h(z) = (k-m) \log |z| + \log |g(z)|$. But the choice of
$\alpha$ in the definition of $h$ implies $k = m$, thereby
proving our claim. Moreover from
$|z|^\alpha = e^{u(z)-v(z)-h(z)} \in L^2(D)$ we
conclude that
$$
u(z) = \alpha \log |z| + \omega(z) \quad \mbox{ where }
\alpha > -1, \,\omega \in C^0 \cap W^{1,2}(D).
$$
Next we perform a blowup to show that $\alpha = m$. For any
$z_0 \in \C \backslash \{0\}$ and $0 < \lambda < 1/|z_0|$ we let
$$
f_\lambda:D_{\frac{1}{\lambda}}(0) \to \R^n,\,
f_\lambda(z) = \frac{1}{\lambda^{\alpha +1}}
\big(f(\lambda z) - f(\lambda z_0)\big).
$$
The $f_\lambda$ have induced metric
$(g_\lambda)_{ij} = e^{2u_\lambda} \delta_{ij}$, where
$$
u_\lambda(z) = u(\lambda z) - \alpha \log \lambda
= \alpha \log |z| + \omega(\lambda z).
$$
Putting $\omega_0 = \omega(0)$ we have
$$
u_\lambda(z) \to \alpha \log |z| + \omega_0 \quad \mbox{ in }
C^0_{loc} \cap W^{1,2}_{loc}(\C \backslash \{0\}).
$$
Furthermore, the Gau{\ss} map of $f_\lambda$ is given by
$G_\lambda(z) = G(\lambda z)$, in particular $DG_\lambda \to 0$
in $L^2_{loc}(\C \backslash \{0\})$. Using the formula
$$
|D^2 f_\lambda|^2 = 2 e^{2u_\lambda} \big(|DG_\lambda|^2 + 2 |Du_\lambda|^2\big),
$$
we obtain by Vitali's theorem
$$
|D^2 f_\lambda|(z) \to \frac{2 e^{\omega_0} \alpha}{|z|^{1-\alpha}}
\quad \mbox{ in } L^2_{loc}(\C \backslash \{0\}).
$$
As $f_\lambda(z_0) = 0$, we can find a sequence $\lambda_k \searrow 0$
such that the $f_{\lambda_k}$ converge in $C^0_{loc}(\C \backslash \{0\})$
and weakly in $W^{2,2}_{loc}(\C \backslash \{0\})$ to a limit map
$f_0:\C \backslash \{0\} \to \R^n$ satisfying $f_0(z_0) = 0$.
After passing to a further subsequence, we can also assume that
$G_{\lambda_k} \to L$ in $W^{1,2}_{loc}(\C \backslash \{0\})$ where
$L \in G(n,2)$ is a constant. It is then easy to see that $f_0$
maps into the plane $L$. Further we have
$$
\langle \partial_i f_0(z),\partial_j f_0(z) \rangle =
e^{2\omega_0}  |z|^{2\alpha} \delta_{ij}.
$$
Using that $f_0$ is locally in $W^{2,2} \cap W^{1,\infty}$ we verify
the identity
$$
\langle \Delta f_0,\partial_j f_0 \rangle =
\partial_i \langle \partial_i f_0,\partial_j f_0 \rangle -
\frac{1}{2} \partial_j \langle \partial_i f_0,\partial_i f_0 \rangle.
$$
Since $f_0$ is conformal, maps into $L$ and has rank two almost everywhere
we see that $f_0$ is harmonic on $\C \backslash \{0\}$. Identifying
$L \cong \C$ by choosing an orthonormal frame $e_{1,2}$, the
conformality relations are
$$
4\, \frac{\partial f_0}{\partial z}
\Big(\overline{\frac{\partial f_0}{\partial \overline{z}}}\Big) =
\left|\frac{\partial f_0}{\partial x}\right|^2 - \left|\frac{\partial f_0}{\partial y}\right|^2
- 2i \left\langle \frac{\partial f_0}{\partial x}, \frac{\partial f_0}{\partial y} \right \rangle
= 0.
$$
Since the two factors on the left are holomorphic, the identity
principle implies that $f_0$ is holomorphic on $\C \backslash \{0\}$,
after replacing $e_1,e_2$ by $e_1,-e_2$ if necessary. Now
$|f'_0(z)| = e^{\omega_0} |z|^\alpha$ and thus for some $\beta \in [0,2\pi)$
$$
f'_0(z) = e^{\omega_0 + i \beta} z^\alpha \quad \mbox{ on }\C \backslash [0,\infty).
$$
As $f'_0$ is single-valued, we must have $\alpha = m \in \N_0$ and
$$
f_0(z) = \frac{e^{\omega_0 + i \beta}}{m+1} \big(z^{m+1}-z_0^{m+1}\big).
$$
In particular, we have the desired expansion $u(z) = m \log |z| + v(z) + h(z)$,
and $u$ satisfies the stated differential equation. Furthermore
$$
|D^2 f|^2 = 2 e^{2u} \big(|DG|^2 + 2 |Du|^2\big) \in L^1(D),
$$
thus $f \in W^ {2,2}(D,\R^n)$. Assuming without loss of generality $f(0) = 0$,
we claim that
$$
\lim_{z \to 0} \frac{|f(z)|}{|z|^{m+1}} = \frac{e^{\omega_0}}{m+1}.
$$
Since $|Df(z)| = |z|^m e^{\omega(z)}$ with $\omega$ bounded, we have
$|f(z)| \leq C |z|^{m+1}$. Now let $z_k \to 0$ be a given sequence.
We can assume that $\zeta_k: = \frac{z_k}{|z_k|} \to \zeta$ with
$|\zeta| = 1$, and compute
\begin{eqnarray*}
\left|\frac{|f(z_k)|}{|z_k|^{m+1}} - \frac{e^{\omega_0}}{m+1}\right|
& = &
\left| \Big|f_{\lambda_k}(\zeta_k) +  \frac{1}{\lambda_k^{m+1}} f(\lambda_k z_0)\Big|
- \Big|\frac{e^{\omega_0 + i \beta}}{m+1} (\zeta^{m+1}-z_0^{m+1})
+ \frac{e^{\omega_0 + i \beta}}{m+1} z_0^{m+1} \Big| \right|\\
& \leq &
\left|f_{\lambda_k}(\zeta_k) - \frac{e^{\omega_0 + i \beta}}{m+1} (\zeta^{m+1}-z_0^{m+1})\right|
+ C |z_0|^{m+1}.
\end{eqnarray*}
Letting $k \to \infty$ we obtain, for a constant $C < \infty$ depending only on
$m$ and $\omega$,
$$
\liminf_{k \to \infty} \left|\frac{|f(z_k)|}{|z_k|^{m+1}} - \frac{e^{\omega_0}}{m+1}\right|
\leq C |z_0|^{m+1}.
$$
This proves our claim since $z_0 \in \C \backslash \{0\}$ was arbitrary. Now
$$
\lim_{\varrho \searrow 0} \frac{\mu_g\big(D_{\varrho}(0)\big)}{\pi r(\varrho)^2} = m+1
\quad \mbox{ where } r(\varrho) = \frac{e^{\omega_0}}{m+1} \varrho^{m+1}.
$$
Choose $\sigma \in (0,1)$ such that $f(z) \neq 0$ for
$z \in \overline{D_\sigma(0)} \backslash \{0\}$, and let $\varrho_{1,2} > 0$
be such that
$$
 \frac{1}{\gamma} r(\varrho_1) = r = \gamma r(\varrho_2), \quad
\mbox{ where } \gamma \in (0,1).
$$
Then for $r > 0$ sufficiently small we have the inclusions
$$
D_{\varrho_1}(0) \subset \big(f^{-1}(B_r(0)) \cap D_\sigma(0)\big) \subset D_{\varrho_2}(0).
$$
It follows that
$$
\gamma^2  \frac{\mu_g(D_{\varrho_1}(0))}{\pi r(\varrho_1)^2} \leq
\frac{f(\mu_g \llcorner D_\sigma(0))(B_r(0))}{\pi r^2} \leq
\frac{1}{\gamma^2} \frac{\mu_g(D_{\varrho_2}(0))}{\pi r(\varrho_2)^2}.
$$
Letting $r \searrow 0$, $\gamma \nearrow 1$ proves that
$\theta^2(f(\mu_g \llcorner D_\sigma(0)),0) = m+1$.
\endproof

A map $f:\Sigma \to \R^n$ is called a branched conformal immersion
(with locally square integrable second fundamental form), if
$f \in W^{2,2}_{{\rm conf},loc}(\Sigma \backslash \mathcal{S},\R^n)$
for some discrete set $\mathcal{S} \subset  \Sigma$ and
$$
\int_{\Omega} |A|^2\,d\mu_g < \infty \quad \mbox{ and } \quad
\mu_g(\Omega) < \infty \quad \mbox{ for all } \Omega \subset \! \subset \Sigma.
$$
The number $m(p)$ as in Theorem \ref{removal} is the
branching order, and $m(p)+1$ is the multiplicity at $p \in \Sigma$.
The map $f$ is unbranched at $p$ if and only if $m(p) = 0$.
For a closed Riemann surface $\Sigma$ and a branched conformal
immersion $f:\Sigma \to \R^n$, consider now
$$
\hat{f} = I_{x_0} \circ f: \Sigma \backslash f^{-1}\{x_0\} \to \R^n,\,
\quad \mbox{ where }
I_{x_0}(x) = x_0 + \frac{x-x_0}{|x-x_0|^2}.
$$
Then $\hat{g} = e^{2v} g$ where $v = - \log |f-x_0|^2$.
The weak Liouville equation says that
$$
\int_{\Sigma} \varphi K_{\hat{g}}\, d\mu_{\hat{g}}
- \int_{\Sigma} \varphi K_g\,d\mu_g =
- \int_{\Sigma} \langle D \log |f-x_0|^2,D\varphi \rangle_g\,d\mu_g \quad
\mbox{ for all } \varphi \in C^\infty_0(\Sigma \backslash f^{-1}\{x_0\}).
$$
A simple computation shows $ |\hat{A}^\circ|^2 d\mu_{\hat{g}} = |A^\circ|^2\,d\mu_g$,
hence by the  Gau{\ss} equation
$$
\frac{1}{4} \int_{\Sigma} \varphi |\hat{H}|^2\,d\mu_{\hat{g}}
- \frac{1}{4} \int_{\Sigma} \varphi |H|^2\,d\mu_g =
- \int_{\Sigma} \langle D \log |f-x_0|^2,D\varphi \rangle_g\,d\mu_g.
$$
Approximating $\varphi \equiv 1$ by cutting off at suitable
radii near each point $p \in f^{-1}\{x_0\}$, we conclude from
the asymptotic information of Theorem \ref{removal}
\begin{equation}
\label{inversion}
{\cal W}(\hat{f}) = {\cal W}(f)
-4\pi \sum_{p \in f^{-1}\{x_0\}} \big(m(p)+1\big).
\end{equation}

\section{Weak compactness of conformal immersions}


\begin{pro}\label{fixcc} Let $\Sigma$ be a closed Riemann surface and
$f_k\in W^{2,2}_{{\rm conf}}(\Sigma,\R^n)$ be a sequence of conformal
immersions satisfying
$$
{\cal W}(f_k) \leq \Lambda < \infty.
$$
Then for a subsequence there exist M\"obius transformations $\sigma_k$ and a
finite set $\mathcal{S} \subset \Sigma$, such that
$$
\sigma_k \circ f_k \to f \quad \mbox{ weakly in }
W^{2,2}_{loc}(\Sigma \backslash \mathcal{S},\R^n),
$$
where $f:\Sigma \to \R^n$ is a branched conformal immersion with
square integrable second fundamental form. Moreover, if $\Lambda < 8\pi$
then $f$ is unbranched and topologically embedded.
\end{pro}

We will use the following standard estimate.

\begin{lem}\label{grad}
Let $\Sigma$ be a two-dimensional, closed manifold with smooth Riemannian
metric $g_0$, and suppose that $u\in W^{1,2}(\Sigma)$
is a weak solution of the equation
$$
-\Delta_{g_0} u = F, \quad \mbox{ where }F \in L^1(\Sigma).
$$
Then for any Riemannian ball $B_r(p)$ and $q\in [1,2)$ we have
$$
\|Du\|_{L^q(B^g_r(p))} \leq C\, r^{\frac{2}{q}-1} \|F\|_{L^1(\Sigma)}
\quad \mbox{ where }C = C(\Sigma,g_0,q) < \infty.
$$
\end{lem}

\proof We may assume that $\|F\|_{L^1(\Sigma)} = 1$ and
$\int_{\Sigma} u\,d\mu_{g_0} = 0$. The function $u$ is given by
$$
u(x) = \int_{\Sigma} G(x,y) F(y)\,d\mu_{g_0}(y),
$$
where $G(x,y)$ is the Riemannian Green function, see Theorem 4.13 in \cite{Aub}.
In particular $G(x,y) = G(y,x)$, and we have the estimate
$$
|D_x G(x,y)| \leq \frac{C}{d(x,y)} \quad \mbox{ where }
C = C(\Sigma,g) < \infty.
$$
By Jensen's inequality we get
\begin{eqnarray*}
\int_{B_r(p)} |Du|^q\,d\mu_{g_0} & \leq &
\int_{B_r(p)} \Big(\int_{\Sigma} |D_x G(x,y)| |F(y)|\,d\mu_{g_0}(y)\Big)^q \,d\mu_{g_0}(x)\\
& \leq & \int_{\Sigma} |F(y)| \int_{B_r(p)}  |D_x G(x,y)|^q\,d\mu_{g_0}(x)\,d\mu_{g_0}(y)\\
& \leq & C \int_{\Sigma} |F(y)| \int_{B_r(p)}  \frac{1}{d(x,y)^q}\,d\mu_{g_0}(x)\,d\mu_{g_0}(y).
\end{eqnarray*}
Now if $d(p,y) < 2r$ we can estimate
$$
\int_{B_r(p)} \frac{1}{d(x,y)^q}\,d\mu_{g_0}(x) \leq
\int_{B_{3r}(y)} \frac{1}{d(x,y)^q}\,d\mu_{g_0}(x) \leq C r^{2-q}.
$$
In the other case $d(p,y) \geq 2r$ we have $d(x,y) \geq r$ on $B_r(p)$,
which implies
$$
\int_{B_r(p)} \frac{1}{d(x,y)^q}\,d\mu_{g_0}(x) \leq
\frac{C}{r^q}\, \mu_g(B_r(p)) \leq C r^{2-q}.
$$
The statement of the lemma follows.
\endproof

{\it Proof of Proposition \ref{fixcc}: }
We may assume $\mu_{g_k} \llcorner |A_k|^2 \to \alpha$ as Radon measures,
and put
$$
\mathcal{S} = \{p\in\Sigma: \alpha(\{p\}) \geq \gamma_n\}.
$$
Choose a smooth, conformal background metric $g_0$
and write $g_k = e^{2u_k} g_0$. Then
$$
\int_{\Sigma} |K_{g_k} e^{2u_k}|\,d\mu_{g_0} =
\int_{\Sigma}|K_{g_k}|\,d\mu_{g_k} \leq
\frac{1}{2}\int_{\Sigma}|A_k|^2\,d\mu_{g_k} \leq C(\Lambda).
$$
From the equation $-\Delta_{g_0} u_k =K_{g_k} e^{2u_k} - K_{g_0}$,
we thus obtain using Lemma \ref{grad} for arbitrary $q \in (1,2)$
the bound
$$
\int_{\Sigma}|D u_k|^q\,d\mu_{g_0} \leq C = C(\Lambda,\Sigma,g_0,q).
$$
By dilating the $f_k$ appropriately we can arrange that
$$
\int_{\Sigma} u_k\,d\mu_{g_0} = 0,
$$
and then get by the Poincar\'{e} inequality, see Theorem 2.34 in \cite{Aub},
$$
\|u_k\|_{W^{1,q}(\Sigma)} \leq C.
$$
In particular, we can assume that $u_k \to u$ weakly in $W^{1,q}(\Sigma)$.
For any $p \notin \mathcal{S}$, we choose conformal coordinates
on a neighborhood $U_\delta(p) \cong D_\delta(0)$, where
$U_\delta(p) \subset \! \subset \Sigma \backslash \mathcal{S}$.
Putting $(g_k)_{ij} = e^{2v_k} \delta_{ij}$ we have
$(g_0)_{ij} = e^{2(v_k-u_k)} \delta_{ij}$ and hence, for
a constant depending on $U_\delta(p)$,
$$
\|v_k\|_{W^{1,q}(U_\delta(p))} \leq C.
$$
Passing to a smaller $\delta > 0$ if necessary, we obtain
from Corollary \ref{key2} the estimate
$$
\|v_k\|_{L^\infty(U_\delta(p))} +
\|D v_k\|_{L^2(U_\delta(p))} \leq C.
$$
Hence we can assume that $v_k$ converges to $v$ on $U_\delta(p)$
weakly in $W^{1,2}$ and pointwise almost everywhere.
But now $|Df_k| = e^{v_k}$ and $\Delta f_k = e^{2v_k}H_k$, where
by assumption
$$
\int_{U_\delta(p)}|H_k|^2e^{2v_k}\,dx dy \leq \Lambda.
$$
Translating the $f_k$ such that $f_k(p) = 0$ for some fixed
$p \in \Sigma \backslash  \mathcal{S}$, we finally obtain
$$
\|f_k\|_{W^{2,2}(\Omega)} \leq C \quad \mbox{ for any }
\Omega \subset \! \subset \Sigma \backslash \mathcal{S}.
$$
In particular the $f_k$ converge weakly in
$W^{2,2}_{loc}(\Sigma \backslash \mathcal{S})$
to some $f \in W^{2,2}_{loc}(\Sigma \backslash \mathcal{S})$,
where $f$ has induced metric $g = e^{2u} g_0$ and
$u \in L^\infty_{loc}(\Sigma \backslash \mathcal{S})$.
If $\limsup_{k \to \infty} \mu_{g_k}(\Sigma) < \infty$,
then $\mu_g(\Sigma) < \infty$ by Fatou's lemma, and the
main statement of Proposition \ref{fixcc} follows from
Theorem \ref{removal}.\\
\\
To prove the statement also in the case $\mu_{g_k}(\Sigma) \to \infty$,
suppose that there is a ball $B_1(x_0)$ with
$f_k(\Sigma) \cap B_1(x_0) = \emptyset$ for all $k$.
Then $\hat{f}_k = I_{x_0} \circ f_k$ converges to $\hat{f} = I_{x_0} \circ f$
weakly in $W^{2,2}_{loc}(\Sigma \backslash \mathcal{S})$,
and $\hat{f}$ has induced metric $\hat{g} = e^{2 \hat{u}}g_0$ where
$\hat{u} = u - \log |f-x_0|^2 \in L^\infty_{loc}(\Sigma \backslash \mathcal{S})$.
Moreover, Lemma 1.1 in \cite{S} yields that
$$
\mu_{\hat{g}_k}(\Sigma)
\leq \Lambda \big(\diam \hat{f}_k(\Sigma)\big)^2 \leq 2 \Lambda.
$$
Thus $\mu_{\hat{g}}(\Sigma) < \infty$ and the result follows as above.
To find the ball $B_1(x_0)$ we employ an argument from \cite{K-S3}.
For $\mu_k = f_k(\mu_{g_k})$ we have by equation (1.3) in \cite{S}
$$
\mu_k(B_R(0)) \leq C R^2 \quad \mbox{ for all } R > 0.
$$
Thus $\mu_k \to \mu$ and $f_k(\mu_{g_k} \llcorner |H_k|^2) \to \nu$ as
Radon measures after passing to a subsequence. Equation 1.4
in \cite{S} implies in the limit
$$
\frac{\mu(B_\varrho(x))}{\varrho^2} + \nu(B_\varrho(x)) \geq c > 0
\quad \mbox{ for all } x \in {\rm spt\,}\mu,\,\varrho > 0.
$$
As shown in \cite{S}, p. 310, the set of accumulation points of
the sets $f_k(\Sigma)$ is just ${\rm spt\,}\mu$. For $R > 0$
to be chosen, let $B_2(x_j)$, $j = 1,\ldots,N$, be a maximal
collection of $2$-balls with centers $x_j \in B_R(0)$, hence
$N \geq R^n/4^n$. Now if ${\rm spt\,}\mu \cap B_1(x_j) \neq \emptyset$
for all $j$, then summation of the inequality over the balls yields
$$
cN \leq \sum_{j=1}^N \Big(\mu\big(B_2(x_j)\big) + \nu\big(B_2(x_j)\big)\Big)
\leq C(\Lambda,n) (R^2 +1).
$$
Therefore ${\rm spt\,}\mu \cap B_1(x_j) = \emptyset$ for some
$j$, if $R = R(\Lambda,n)$ is suffciently large. The
additional conclusions in the case $\Lambda < 8\pi$ are
clear from formula (\ref{inversion}) and Theorem \ref{removal}.
\endproof

The following existence result is proved independently in a
recent preprint by Rivi\`{e}re \cite{R}. It
extends previous work of Kuwert and Sch\"atzle \cite{K-S4}.
In their paper, it is shown that the minimizers are actually
smooth.

\begin{cor} Let $\Sigma$ be a closed Riemann surface such that
$$
\beta^n_\Sigma =
\inf \{{\cal W}(f): f \in W^{2,2}_{{\rm conf}}(\Sigma,\R^n)\} < 8\pi.
$$
Then the infimum $\beta^n_\Sigma$ is attained.
\end{cor}

We now generalize Proposition \ref{fixcc} to the case
of varying Riemann surfaces. The following standard
lemma will be useful, see \cite{D-K} for a proof.

\begin{lem}\label{isothermal}
Let $g_k,g$ be smooth Riemannian metrics on a surface $M$,
such that $g_k \to g$ in $C^{s,\alpha}(M)$, where $s \in \N$,
$\alpha \in (0,1)$. Then for each $p \in M$ there exist
neighborhoods $U_k, U$ and smooth conformal diffeomorphisms
$\varphi_k:D \to U_k$, such that $\varphi_k \to \varphi$
in $C^{s+1,\alpha}(\overline{D},M)$.
\end{lem}

%
%

\begin{thm}\label{convergence2}
Let $f_k \in W^{2,2}(\Sigma_k,\R^n)$ be conformal immersions
of compact Riemann surfaces of genus $p$. Assume that the
$\Sigma_k$ converge to $\Sigma$ in moduli space, i.e.
$\phi_k^\ast(\Sigma_k) \to \Sigma$ as complex structures
for suitable diffeomorphisms $\phi_k$, and that
$$
{\cal W}(f_k) \leq \Lambda < \infty.
$$
Then there exist a branched conformal immersion $f:\Sigma \to \R^n$
with square integrable second fundamental form, a finite set
$\mathcal{S} \subset M$ and M\"obius transformations
$\sigma_k$, such that for a subsequence
$$
\sigma_k \circ f_k \circ \phi_k \to f \quad \mbox{ weakly in }
W^{2,2}(\Sigma \backslash \mathcal{S},\R^n).
$$
\end{thm}

The convergence of the complex structures implies
that $\phi_k^\ast g_{0,k} \to g_0$, where $g_{0,k}$, $g$
are the suitably normalized, constant curvature metrics
in $\Sigma_k$, $\Sigma$, see chapter 2.4 in \cite{T}.
The proof is now along the lines of Proposition \ref{fixcc},
using the local conformal charts from Lemma \ref{isothermal}.

%
%
%

\section{The energy of surfaces diverging in moduli space}

\subsection{H\'elein's convergence theorem}

The following result is essentially Theorem 5.1.1 in \cite{H},
except that the constant $8\pi/3$ is replaced here by $\gamma_n$,
exploiting the estimate from \cite{M-S}. At the end of the
subsection we will show that $\gamma_n$ is in fact optimal.

\begin{thm}\label{Helein} Let $f_k\in W^{2,2}_{{\rm conf}}(D,\R^n)$
be a sequence of conformal immersions with induced metrics
$(g_k)_{ij} = e^{2u_k} \delta_{ij}$, and assume
$$
\int_D |A_{f_k}|^2\,d\mu_{g_k} \leq \gamma <
\gamma_n =
\begin{cases}
8\pi & \mbox{ for } n = 3,\\
4\pi & \mbox{ for }n \geq 4.
\end{cases}
$$
Assume also that $\mu_{g_k}(D) \leq C$ and $f_k(0) = 0$.
Then $f_k$ is bounded in $W^{2,2}_{loc}(D,\R^n)$, and there
is a subsequence such that one of the following two alternatives
holds:
\begin{itemize}
\item[{\rm (a)}] $u_k$ is bounded in $L^\infty_{loc}(D)$ and
$f_k$ converges weakly in $W^{2,2}_{loc}(D,\R^n)$ to a conformal
immersion $f \in W^{2,2}_{{\rm conf},loc}(D,\R^n)$.
\item[{\rm (b)}] $u_k \to - \infty$ and $f_k \to 0$ locally uniformly on $D$.
\end{itemize}
\end{thm}

\proof By Corollary \ref{key} there is a solution $v_k$ of the
equation $- \Delta v_k = K_{g_k} e^{2 u_k}$ satisfying
$$
\|v_k\|_{L^\infty(D)} + \|Dv_k\|_{L^2(D)}
\leq C(\gamma) \int_D |A_{f_k}|^2\,d\mu_{g_k}.
$$
Clearly $h_k = u_k - v_k$ is harmonic on $D$. Now
$$
\int_D e^{2u_k^+} =
|\{u_k \leq 0\}| + \int_{\{u_k > 0\}} e^{2u_k} \leq C,
$$
and hence by Jensen's inequality
$$
\int_D u_k^+ \leq C.
$$
For $\dist(z,\partial D) \geq r$ where $r \in (0,1)$ we get
$$
h_k(z) =
\frac{1}{\pi r^2} \int_{D_r(z)} (u_k - v_k)  \leq
\frac{1}{\pi r^2} \int_D u_k^+  + \|v_k\|_{L^\infty(D)}
\leq C(\gamma,r).
$$
Thus $u_k = v_k + h_k$ is locally bounded from above, which
implies that the sequence $f_k$ is bounded in $W^{1,\infty}_{loc}(D,\R^n)$.
As $\Delta f_k = e^{2u_k} H_{f_k}$, we further have for $\Omega = D_{1-r}(0)$
$$
\int_{\Omega} |\Delta f_k|^2\, = \int_{\Omega} e^{2u_k} |H_{f_k}|^2\,d\mu_{g_k} \leq
C(\gamma,r) \int_{\Omega} |A_{f_k}|^2\,d\mu_{g_k} \leq C(\gamma,r).
$$
Thus $f_k$ is also bounded in $W^{2,2}_{loc}(D,\R^n)$ and converges,
after passing to a subsequence, weakly to some
$f \in W^{2,2}_{loc} \cap W^{1,\infty}_{loc}(D,\R^n)$.
Now if $\int_D u_k^- \leq C$, then for $\dist(z,\partial D) \geq r$
$$
h_k(z) =
\frac{1}{\pi r^2} \int_{D_r(z)} (u_k - v_k)  \geq
- \frac{1}{\pi r^2} \int_D u_k^-  - \|v_k\|_{L^\infty(D)}
\geq - C(\gamma,r).
$$
Thus $u_k = v_k + h_k$ is bounded in $L^\infty_{loc} \cap W^{1,2}_{loc}(D)$,
and $u_k$ converges pointwise to a function $u \in L^\infty_{loc}(D)$.
We conclude
$$
g_{ij} = \langle \partial_i f, \partial_j f \rangle = e^{2u} \delta_{ij},
$$
which means that $f$ is a conformal immersion as claimed in case (a).
We will now show that $\int_D u_k^- \to \infty$ implies alternative (b).
Namely, we then have
$$
h_k(0) = \frac{1}{\pi} \int_D (u_k-v_k) \to - \infty.
$$
As $C(\gamma,r) - h_k \geq 0$ on $\Omega$, we get by the Harnack inequality
$$
\sup_{\Omega'} h_k \leq
\frac{1}{C(r)} \inf_{\Omega'} h_k + C(\gamma,r) \to - \infty \quad
\mbox{ where } \Omega' = D_{1-2r}(0).
$$
Thus $u_k = v_k + h_k  \to - \infty$ and $f_k \to 0$ locally uniformly on $D$.
\endproof

Applying Lemma \ref{isothermal}, we get a version of H\'elein's
theorem for conformal immersions with respect to a convergent
sequence of metrics.

\begin{cor}\label{Helein2} The statement of Theorem \ref{Helein}
continues to hold for immersions $f_k \in W^{2,2}(D,\R^n)$
with induced metric $g_k = e^{2u_k} g_{0,k}$, if the
$(g_{0,k})_{ij}$ converge to $\delta_{ij}$ smoothly on $\overline{D}$.
\end{cor}

Relating to Remark 5.1.3 in \cite{H}, we now show that the constant
$4\pi$ in Theorem \ref{Helein} is optimal for $n \geq 4$.
For $\varepsilon > 0$ we consider the conformally immersed
minimal disks
$$
f_\varepsilon:D \to \C^2,\,f_\varepsilon(z) =
\Big(\frac{1}{2} z^2,\varepsilon z \Big).
$$
We compute $(g_\varepsilon)_{ij} = e^{2u_\varepsilon} \delta_{ij}$ where
$u_\varepsilon(z) = \frac{1}{2} \log (|z|^2 + \varepsilon^2)$, and further
$$
\int_D |A_{f_\varepsilon}|^2\,d\mu_{g_\varepsilon} =
- 2 \int_D K_{g_\varepsilon}\,d\mu_{g_\varepsilon} =
2 \int_{\partial D} \frac{\partial u_\varepsilon}{\partial r}\,ds
= \frac{4\pi}{1+\varepsilon^2} < 4\pi.
$$
As $f_\varepsilon(z) \to (\frac{1}{2}z^2,0)$ for $\varepsilon \searrow 0$,
none of the two alternatives (a) or (b) is satisfied. For the
optimality of $\gamma_3 = 8\pi$ we also follow \cite{M-S} and
consider Enneper's minimal surface
$$
f:\C \to \R^3,\,f(z) =
\frac{1}{2}\,{\rm Re\,}\Big(z-\frac{1}{3}z^3,i\big(z+ \frac{1}{3}z^3\big),z^2\Big).
$$
We have
$f_\lambda(z) = \frac{1}{\lambda^3} f(\lambda z) \to - \frac{1}{6} (z^3,0) \in \C \times \R = \R^3$
as $\lambda \nearrow \infty$. Restricting $f_\lambda$ to $D$ yields conformally
immersed disks with $\int_D |A_{f_\lambda}|^2\,d\mu_{g_\lambda} < 8\pi$.

\subsection{The case of tori}

The following was proved in \cite{K-S3} for $n=3$, and for
$n = 4$ with bound $\min(8\pi,\beta^4_1 + \frac{8\pi}{3})$.

\begin{thm}\label{genus1} Let $\Sigma_k$ be complex tori which diverge
in moduli space. Then for any sequence of conformal immersions
$f_k \in W^{2,2}_{{\rm conf}}(\Sigma_k,\R^n)$ we have
$$
\liminf_{k \to \infty}{\cal W}(f_k) \geq 8\pi.
$$
\end{thm}


\proof We may assume that $\Sigma_k =  \C/\Gamma_k$ where
$\Gamma_k = \Z \oplus \Z(a_k+ib_k)$ is normalized by
$0 \leq a_k \leq \frac{1}{2}$, $a_k^2 + b_k^2 \geq 1$ and $b_k > 0$.
We also assume that the $f_k:\Sigma_k \to \R^n$ satisfy
$$
\frac{1}{4} \limsup_{k \to \infty} \int_{\Sigma_k} |A_{f_k}|^2\,d\mu_{g_k}
= \limsup_{k \to \infty} {\cal W}(f_k) \leq \Lambda < \infty.
$$
We lift the $f_k$ to $\Gamma_k$-periodic maps from $\C$ into $\R^n$.
Theorem \ref{removal} shows that $f_k$ is not constant on any circle
$C_v = [0,1] \times \{v\}$, $v \in \R$. Hence by passing to
$\frac{1}{\lambda_k} (f(u,v+v_k) - f(0,v_k))$ for suitable
$\lambda_k > 0$, $v_k \in [0,b_k)$, we may assume that
$$
1 = \diam f_k(C_0) \leq \diam f_k(C_v) \quad \mbox{ for all }v \in \R,
\quad \mbox{ and } \quad f_k(0,0) = 0.
$$
Arguing as in the proof of Proposition \ref{fixcc}, we obtain
$B_1(x_0) \subset \R^n$ such that $f_k(\Sigma_k) \cap B_1(x_0) = \emptyset$
for all $k$. For $\hat{f}_k = I_{x_0} \circ f_k$ we have
$\hat{f}_k(\Sigma_k) \subset \overline{B_1(x_0)}$, and Lemma 1.1 in
\cite{S} implies an area bound $\mu_{\hat{g}_k}(\Sigma_k) \leq C$.
Up to a subsequence, we have
$\mu_{\hat{g}_k} \llcorner |A_{\hat{f}_k}|^2 \to \alpha$
as Radon measures on the cylinder $C = [0,1] \times \R$.
The set $\mathcal{S} = \{w \in C: \alpha(\{w\}) \geq \gamma_n\}$
is discrete, and
$$
\varrho(w) = \inf \{\varrho > 0: \alpha(D_\varrho(w)) \geq \gamma_n\} > 0
\quad \mbox{ for } w \in \Omega = C \backslash \mathcal{S}.
$$
Now $\hat{f}_k$ converges locally uniformly in $\Omega$ either
to a conformal immersion, or to a point $x_1 \in \R^n$.
This follows from Theorem \ref{Helein} together with a continuation
argument, using that $\varrho(w)$ is lower semicontinuous and
hence locally bounded from below. Note
$$
\hat{f_k}(C_0) \subset I_{x_0}\big(\overline{B_1(0)}\big) \subset
\R^n \backslash B_\theta(x_0) \quad
\mbox{ where }\theta = \frac{1}{|x_0|+1} > 0.
$$
In the second alternative we thus get $|x_1 -x_0| \geq \theta > 0$,
and $f_k|_{C_v}$ converges uniformly to the point $I_{x_0}(x_1)$
for any $C_v \subset \Omega$, in contradiction to $\diam f_k(C_v) \geq 1$.
Therefore $\hat{f}_k$ converges to a conformal immersion
$\hat{f}:\Omega \to \R^n$. Now the assumption that $\Sigma_k$
diverges in moduli space yields that $b_k \to \infty$, so that
$\hat{f}:\Omega \to \R^n$ has second fundamental form in $L^2(C)$
and finite area. Applying Theorem \ref{removal} to the points
at $v = \pm \infty$ we see that $f(C_v) \to x_{\pm} \in \R^n$
for $v \to \pm \infty$. Let us assume that $x_{+} \neq x_0$.
Then for any $\varepsilon > 0$ we find a $\delta > 0$ with
$I(B_\delta(x_{+})) \subset B_\varepsilon(I(x_+))$.
Choosing $v < \infty$ large such that
$\hat{f}(C_v) \subset B_{\frac{\delta}{2}}(x_{+})$,
we get for sufficiently large $k$
$$
f_k(C_v) = I\big(\hat{f}_k(C_v)\big) \subset
I\big(B_\delta(x_{+})\big) \subset  B_\varepsilon(I(x_+)).
$$
Taking $\varepsilon = \frac{1}{3}$ yields a contradiction to
$\diam f_k(C_v) \geq 1$. This shows $x_\pm = x_0$ and in particular
$\theta^2(\hat{\mu},x_0) \geq 2$ where $\hat{\mu} = \hat{f}(\mu_{\hat{g}})$.
We conclude from the Li-Yau inequality, see Section 2.2,
$$
\liminf_{k \to \infty} {\cal W}(f_k) =
\liminf_{k \to \infty} {\cal W}(\hat{f}_k) \geq
{\cal W}(\hat{f}) \geq 8\pi.
$$
\endproof

\subsection{The case of genus $p \geq 2$}
We first collect some facts about degenerating Riemann
surfaces from \cite{B,Hum}. By definition, a compact Riemann
surface with nodes is a compact, connected Hausdorff space
$\Sigma$ together with a finite subset $N$, such that
$\Sigma \backslash N$ is locally  homeomorphic to $D$,
while each $a \in N$ has a neighborhood  homeomorphic to
$\{(z,w) \in \C^2: zw = 0,\, |z|, |w| < 1\}$.
Moreover, all transition functions are required to
be holomorphic. The points in $N$ are called nodes. Each
component $\Sigma^i$ of $\Sigma \backslash N$ is contained
in a compact Riemann surface $\overline{\Sigma}^i$, which
is given by adding points to the punctured coordinate disks
at the nodes. We have $q \leq \nu +1$, where $q$ and $\nu$
are the number of the components and the nodes, respectively.
We denote by $p_i$ the genus of $\overline{\Sigma}^i$ and $\nu_i$
the number of punctures of $\Sigma^i$. If $2p_i + \nu_i \geq 3$
or equivalently $\chi(\Sigma^i) < 0$, then $\Sigma^i$ carries
a unique conformal, complete metric having constant
curvature $-1$. With respect to this metric, the surface
has cusps at the punctures and area $4\pi(p_i-1+\nu_i)$.\\
\\
Next let $\Sigma_k$ be a sequence of compact Riemann surfaces of
fixed genus $p \geq 2$, with hyperbolic metrics $h_k$.
By Proposition 5.1 in \cite{Hum}, there exists a compact Riemann surface
$\Sigma$ with nodes $N = \{a_1,\ldots,a_r\}$, and for each $k$
a maximal collection $\Gamma_k = \{\gamma_k^1,\ldots,\gamma_k^r\}$
of pairwise disjoint, simply closed geodesics in $\Sigma_k$
with $\ell^j_k = L(\gamma_k^j) \to 0$, such that after passing
to a subsequence the following holds:
\begin{itemize}
\item[{\rm (1)}] $p - \sum_{i=1}^q p_i = \nu+1-q \geq 0$.
\item[{\rm (2)}] There are maps $\varphi_k \in C^0(\Sigma_k,\Sigma)$,
such that $\varphi_k: \Sigma_k \backslash \Gamma_k \to \Sigma \backslash N$
is diffeomorphic and $\varphi_k(\gamma_k^j) = a_j$ for $j = 1,\ldots,r$.
\item[{\rm (3)}] For the inverse diffeomorphisms
$\psi_k:\Sigma \backslash N \to \Sigma_k \backslash \Gamma_k$,
we have $\psi_k^\ast h_k \to h$ in $C^\infty_{loc}(\Sigma \backslash N)$.
\end{itemize}

In the following we consider a sequence of conformal immersions
$f_k \in W^{2,2}(\Sigma_k,\R^n)$ with ${\cal W}(f_k) \leq \Lambda$,
and we assume that the hyperbolic surfaces $(\Sigma_k,h_k)$
converge to a surface with nodes $(\Sigma,N)$ as described above.

\begin{lem}\label{genusp} There exist branched conformal
immersions $f^i:\overline{\Sigma}^i \to \R^n$, finite sets
$\mathcal{S}_i \subset \Sigma^i$ and M\"obius transformations
$\sigma^i_k$, such that
$$
\sigma^i_k \circ f_k \circ \psi_k|_{\Sigma^i}  \to f^i
\quad \mbox{ weakly in } W^{2,2}_{loc}(\Sigma^i \backslash \mathcal{S}_i,\R^n)
\mbox{ for }i = 1,\ldots,q.
$$
Replacing $f_k$ by $\sigma_k \circ f_k$ for suitable M\"obius transformations
$\sigma_k$ we can take $\sigma^1_k = {\rm id}$ and
$$
\sigma^i_k(y) = I_{x_i}\Big(\frac{y-y^i_k}{\lambda^i_k}\Big)
\quad \mbox{ where }x_i \in \R^n,\,y^i_k = (f_k \circ \psi_k)(b_i)
\mbox{ for }b_i \in \Sigma^i \mbox{ and }\lambda^i_k > 0,
$$
for $i = 2,\ldots,q$. Further the maps $\sigma^i_k \circ f_k$ are
uniformly bounded, and ${\cal W}(f^i) \geq \beta^n_{p_i}$.
\end{lem}

\proof By the Gau{\ss}-Bonnet formula, the second fundamental form is bounded by
$$
\int_{\Sigma_k} |A_{f_k}|^2\,d\mu_{f_k} \leq C(\Lambda,p) < \infty.
$$
The maps $f_k \circ \psi_k:\Sigma \backslash N \to \R^n$ are
conformal immersions with respect to the metric
$\psi_k^\ast h_k$, which converges to $h$ in
$C^\infty_{loc}(\Sigma \backslash N)$. Let $\xi_i \subset \Sigma^i$
be an embedded arc, which is subdivided into
$\xi_i^1,\ldots,\xi_i^m$. We can choose a subsequence
and $j_0 \in \{1,\ldots,m\}$ with
$$
\diam (f_k \circ \psi_k)(\xi_i^{j_0}) =
\min_{1 \leq j \leq m}  \diam (f_k \circ \psi_k) (\xi_i^j)
=: \lambda^i_k.
$$
We have $\lambda^i_k > 0$ by Theorem \ref{removal}. Select
$b_i \in \xi_i^{j_0}$, and define the maps
$$
f^i_k:\Sigma_k \to \R^n,\,f^i_k(p) = \frac{f_k(p)-y^i_k}{\lambda^i_k}
\quad \mbox{ where } y^i_k = (f_k \circ \psi_k)(b_i).
$$
As in Proposition \ref{fixcc}, we can choose
$B_1(x_i) \subset \R^n$ with $f^i_k(\Sigma_k) \cap B_1(x_i) = \emptyset$
for all $k$. Applying \eqref{diameter} to $I_{x_i} \circ f^i_k$ yields
$$
\mu_{I_{x_i} \circ f^i_k}(\Sigma_k) \leq C < \infty.
$$
Now consider the maps
$$
\hat{f}^i_k = I_{x_i} \circ f^i_k \circ \psi_k|_{\Sigma^i}:\Sigma^i \to \R^n.
$$
We can assume that $\mu_{\hat{f}^i_k} \llcorner |A_{\hat{f}^i_k}|^2$
converges to $\alpha$ as Radon measures, and put
$$
\mathcal{S}_i = \{p \in \Sigma^i: \alpha(\{p\}) \geq \gamma_n\}.
$$
Corollary \ref{Helein2} implies that, away from $\mathcal{S}_i$, the
$\hat{f}^i_k$ subconverge locally uniformly either to a conformal immersion,
or to a point $x_1 \in \R^n$. As in Theorem \ref{genus1}
$$
\hat{f}^i_k (\xi_i^{j_0}) \subset I_{x_i} (\overline{B_1(0)}) \subset
\overline{B_1(x_i)} \backslash B_{\theta_i}(x_i) \quad \mbox{ where }
\theta_i = \frac{1}{|x_i|+1} > 0.
$$
Therefore in the second alternative we get $|x_1-x_i| \geq \theta_i$,
and $f^i_k \circ \psi_k$ converges to $I_{x_i}(x_1)$ locally
uniformly on $\Sigma^i \backslash \mathcal{S}_i$.
But for $m > \frac{C(\Lambda,p)}{\gamma_n}$ there is a $j \in \{1,\ldots,m\}$
with $\xi_i^j \cap \mathcal{S}_i = \emptyset$, and we conclude
$1 \leq \diam (f^i_k \circ \psi_k)(\xi_i^j) \to 0$, a contradiction.
Therefore $\hat{f}^i_k$ converges locally
uniformly and weakly in $W^{2,2}_{loc}(\Sigma^i \backslash \mathcal{S}_i,\R^n)$
to $f^i \in W^{2,2}_{{\rm conf},loc}((\Sigma^i \backslash \mathcal{S}_i,\R^n)$.
Furthermore, Theorem \ref{removal} shows that $f^i$ extends as a branched
conformal immersion to $\overline{\Sigma}^i$. Applying the argument
for $i = 2,\ldots,q$ with $f_k$ replaced by $\sigma^1_k \circ f_k$
yields the second statement of the lemma. Finally, the inequality
${\cal W}(f^i) \geq \beta^n_{p_i}$ is clear when $f^i$ is unbranched,
otherwise we get ${\cal W}(f^i) \geq 8\pi > \beta^n_{p_i}$
from the Li-Yau inequality (\ref{inversion}) in connection with
\cite{K}.
\endproof

For our last result we need more details on degenerating hyperbolic
surfaces. For $\ell > 0$ we define a reference cylinder
$C(\ell) = [0,1] \times [-T(\ell),T(\ell)]$ with metric $g_\ell$, where
$$
T(\ell) = \frac{1}{\ell}\, \arccot \big(\sinh \frac{\ell}{2}\big)
\quad \mbox{ and } \quad
g_\ell(s,t) = \frac{\ell^2}{\cos^2 \ell t}\,(ds^2 + dt^2).
$$
The map $(s,t) \mapsto i e^{\ell(s+it)}$ yields an isometry between
$(C(\ell),g_\ell)$ and the sector in the upper halfplane given
by $1 \leq r \leq e^\ell$,
$|\theta - \frac{\pi}{2}| \leq \arccot \big(\sinh \frac{\ell}{2}\big)$.
The circles $c_t = \{(s,t): s \in [0,1]\}$ have constant geodesic
curvature $\varkappa_{g_\ell}(t) = \sin \ell t$ and length
$L_{g_\ell}(t) = \ell/\cos \ell t$. We note
$$
\lim_{\ell \searrow 0} \varkappa_{g_\ell}\big(\pm(T(\ell)-t)\big) = 1
\quad \mbox{ and }\quad
\lim_{\ell \searrow 0} L_{g_\ell}\big(\pm(T(\ell)-t)\big) = \frac{1}{t + \frac{1}{2}}
\quad \mbox{ for any } t > 0.
$$
Now let $\gamma_k \subset \Sigma_k$ be a sequence of geodesics
with length $\ell_k \to 0$, corresponding to the node $a \in \Sigma$.
By the collar lemma, see \cite{Hum}, there is an isometric embedding
$$
\big(C(\ell_k),g_{\ell_k}\big) \subset (\Sigma_k,h_k),
$$
with $c_0$ corresponding to $\gamma_k$. Clearly $T_k = T(\ell_k) \to \infty$.
We will need the following property of the
construction in \cite{Hum}: for any $t \in [0,\infty)$ there is a compact
set $K_t \subset \Sigma \backslash N$ such that
\begin{equation}
\label{hummelcondition}
\varphi_k([0,1] \times [T_k-t,T_k]) \subset K_t \mbox{ for all }k \in \N.
\end{equation}
For this we refer to Section 4 in \cite{Z}.


\begin{thm} Let $\Sigma_k$ be sequence of compact Riemann surfaces
of genus $p \geq 2$, which diverges in moduli space.
Then for any sequence of conformal
immersions $f_k \in W^{2,2}_{{\rm conf}}(\Sigma_k,\R^n)$ we have
$$
\liminf_{k \to \infty} {\cal W}(f_k) \geq \min(8\pi,\omega^n_p).
$$
\end{thm}

\proof We first consider the case $q = \nu +1$, hence
$p = p_1 + \ldots + p_q$. By Lemma \ref{genusp} we have,
away from a finite set of points, $f_k \circ \psi_k \to f^1$
weakly on $\Sigma^1$ and
$$
\frac{f_k \circ \psi_k - y^i_k}{\lambda^i_k} \to I_{x_i} \circ f^i
\quad \mbox{ weakly on } \Sigma^i \mbox{ for } i = 2,\ldots,q.
$$
Now if $f^i$ attains $x_i$ with multiplicity two or
more, then the Li-Yau inequality (\ref{inversion}) yields
$$
\liminf_{k \to \infty} {\cal W}(f_k) \geq {\cal W}(f^i) \geq  8\pi,
$$
Otherwise we obtain, again by (\ref{inversion}),
$$
\lim_{k \to \infty} {\cal W}(f_k) \geq
{\cal W}(f^1) + \sum_{i = 2}^q {\cal W}(I_{x_i} \circ f^i)
\geq \beta^n_{p_1} + \sum_{i = 2}^q (\beta^n_{p_i} - 4\pi) \geq \omega^n_p.
$$
In the case $q < \nu +1$ there must be a node which does not
disconnect $\Sigma$. After renumbering we can chose components
$\Sigma^1,\ldots,\Sigma^s$, and for each $\Sigma^i$ two punctures
$a_i^{\pm}$ such that $a_i^{+},a_{i+1}^{-}$ correspond to the same
node $a_i$ for $i = 1,\ldots,s$; here $a_{s+1}^{-} = a_1^{-}$.
We say that a puncture $a_i^{\pm}$ is good, if either $i=1$
or $f^i(a_i^{\pm}) \neq x_i$. If both $a_i^{+}$ and $a_i^{-}$
are not good, then the theorem follows with lower bound $8\pi$ by
the Li-Yau inequality (\ref{inversion}). Therefore, omitting
subscripts we can assume that there is a node $a$ at which
both punctures are good.\\
\\
Using the collar embedding we now choose $\tau_k \in [-T_k,T_k]$
with
$$
\diam f_k(c_{\tau_k}) =
\min_{t \in [-T_k,T_k]} \diam f_k(c_{t}) =: \delta_k.
$$
The result follows as in Theorem \ref{genus1} once we can show
that for a subsequence
\begin{equation}
\label{finalclaim}
\lim_{k \to \infty} |T_k \pm \tau_k| = \infty.
\end{equation}
For fixed $t \in [0,\infty)$ the curves $\varphi_k(c_{T_k-t})$
are contained in the compact set $K_t \subset \Sigma \backslash N$.
Since $\psi_k^\ast h_k$ converges to $h$ smoothly on $K_t$,
we can assume that the curves converge smoothly to a limiting
curve $\beta_t$ in $K_t$ with length $L_h(\beta_t) = (t+\frac{1}{2})^{-1}$.
Now if $\varphi_k(c_{T_k-t}) \subset \Sigma^1$ we have
$$
\diam f_k(c_{T_k-t}) = \diam (f_k \circ \psi_k)(\varphi_k(c_{T_k-t})) \to
\diam f^1(\beta_t).
$$
By Theorem \ref{removal}, we see $\diam f^1(\beta_t) > 0$ for any
$t \in [0,\infty)$. On the other hand
$$
\limsup_{k \to \infty} \delta_k \leq \limsup_{k \to \infty} \big(\diam f_k(c_{T_k-t})\big)
= \diam f^1(\beta_t).
$$
Letting $t \to \infty$ we conclude $\lim_{k \to \infty} \delta_k = 0$
by continuity of $f^1$, which proves claim (\ref{finalclaim}).
In the remaining case  $\varphi_k(c_{T_k-t}) \subset \Sigma^i$ for some
$i \geq 2$ we compute similarly
$$
\frac{\diam f_k(c_{T_k-t})}{\lambda^i_k} =
\diam (I_{x_i} \circ f_k \circ \psi_k)(\varphi_k(c_{T_k-t})) \to
\diam (I_{x^i} \circ f^i)(\beta_t) > 0,
$$
and further
$$
\limsup_{k \to \infty} \frac{\delta_k}{\lambda^i_k}
\leq \limsup_{k \to \infty} \frac{\diam f_k(c_{T_k-t})}{\lambda^i_k}
= \diam (I_{x^i} \circ f^i)(\beta_t).
$$
Again letting $t \to \infty$ we see that $\delta_k/\lambda^i_k \to 0$,
using the fact that the puncture is good, i.e. $f^i(a) \neq x_i$.
Thus (\ref{finalclaim}) holds also for $i \geq 2$, and the theorem
is proved.
\endproof

The constants $\beta^n_p$ and hence $\omega^n_p$ are not known explicitely.
The Willmore conjecture in $\R^n$ would imply that
$\omega^n_2 = 4\pi(\pi -1) > 8\pi$. The inequality
$\omega^n_p > 8\pi$ holds at least for large $p$, since
$\beta^n_p \to 8\pi$ as $p \to \infty$ by \cite{K-L-S},
as noted in the introduction.


\vspace{1cm}
\begin{tabular}{ll}
{\sc Ernst Kuwert }&{\sc Yuxiang Li}\\
{\sc Mathematisches Institut}&{\sc Department of Mathematical Sciences}\\
{\sc Albert-Ludwigs-Universit\"at Freiburg}&{\sc Tsinghua University}\\
{\sc Eckerstra{\ss}e 1, D-79104 Freiburg}& {\sc Beijing 100084, P.R. China}\\
{\tt ernst.kuwert@math.uni-freiburg.de} & {\tt yxli@math.tsinghua.edu.cn}
\end{tabular}


{\small\begin{thebibliography}{2}

\bibitem[Aub]{Aub} T. Aubin: Nonlinear analysis on
manifolds. Monge-Amp\`{e}re equations. {\em Grundlehren
der mathematischen Wissenschaften} {\bf 252}, Springer
Verlag, New York 1982.

\bibitem[B-K]{B-K}
M. Bauer and E. Kuwert:
Existence of minimizing Willmore surfaces of prescribed genus.
{\em Int. Math. Res. Not.}, {\bf 10} (2003), 553-576.

\bibitem[B]{B} L. Bers: Spaces of degenerating Riemann surfaces.
{\em Discontinuous groups and Riemann surfaces, Ann. of Math.
Studies} {\bf 79}, 43--55, Princeton Univ. Press, Princeton N.J. (1974).

\bibitem[C]{C} B. Y. Chen: Some conformal invariants
of submanifolds and their application, {\em
Bollettino della
Unione Matematica Italiana} {\bf 10} (1974), 380--385.

\bibitem[C-L-M-S]{C-L-M-S} R. Coifman, P.L. Lions, Y. Meyer, S. Semmes:
Compensated compactness and Hardy spaces, {\em J. Math. Pures
Appl. }{\bf 72}, 247--286.

\bibitem[D-K]{D-K}D. DeTurck and J. Kazdan:
Some regularity theorems in Riemannian geometry.
{\em Ann. Sci. šŠcole Norm. Sup. (4)} {\bf 14}
 (1981),  249--260.

\bibitem[H]{H}F. H\'elein: Harmonic maps,
conservation laws and moving frames.
Translated from the 1996 French original.
With a foreword by James Eells. Second edition.
Cambridge Tracts in Mathematics, 150.
Cambridge University Press, Cambridge, 2002.


\bibitem[Hum]{Hum}C. Hummel: Gromov's compactness
theorem for pseudo-holomorphic curves. {\em Progress in
Mathematics} {\bf 151},
Birkh\"auser Verlag, Basel (1997)

\bibitem[K]{K}
R. Kusner: Comparison surfaces for the Willmore problem.
{\em Pacific J. Math.}, {\bf 138} (1989), 317--345.

\bibitem[LY]{LY} P. Li and S.T. Yau: A new conformal
invariant and its applications to the Willmore conjecture
and the first eigenvalue on compact surfaces,
{\em Invent. Math} {\bf 69} (1982), 269--291.

\bibitem[K-L-S]{K-L-S} E. Kuwert, Y. Li and
R. Sch\"atzle: The large genus limit of the
infimum of the Willmore energy. {\em Amer. J.
Math.} {\bf 132} (2010), 37--51.

\bibitem[K-S1]{K-S2} E. Kuwert and R. Sch\"atzle: Removability
of point singularities of Willmore surfaces, {\em Ann. of Math.}
{\bf 160} (2004), 315--357.

\bibitem[K-S2]{K-S3} E. Kuwert and R. Sch\"atzle: Closed
surfaces with bounds on their Willmore energy, {\em Preprint
Centro di Ricerca Matematica Ennio De Giorgi, Pisa} 2008.

\bibitem[K-S3]{K-S4} E. Kuwert and R. Sch\"atzle: Minimizers
of the Willmore functional under fixed conformal class,
{\em Manuscript}, 2008.

\bibitem[L]{L} J. Langer: A compactness theorem for
surfaces with $L^p$-bounded second fundamental form,
{\em Math. Ann. }{\bf 270} (1985), 223--234.


\bibitem[M-S]{M-S} S. M\"uller and V. \v{S}ver\'ak: On surfaces
of finite total curvature, {\em J. Differential Geom.}
{\bf 42} (1995),
229--258.

\bibitem[R]{R} T. Rivi\`{e}re: Variational principles for immersed
surfaces with $L^2$-bounded second fundamental form,
arXiv:1007.2997v1 (2010).

\bibitem[S-U]{S-U} R. Schoen and K. Uhlenbeck: Boundary regularity
and the Dirichlet problem for harmonic maps, {\em J. Differential Geom. }
{\bf 18} (1983), 253--268.

\bibitem[S]{S} L. Simon: Existence of surfaces minimizing
the Willmore functional, {\em Commun. Anal. and Geom.}
{\bf 1} (1993),
281--326.

\bibitem[T]{T} A. Tromba: Teichm\"uller theory in Riemannian geometry,
{\em Lectures in Mathematics ETH Z\"urich}, Birkh\"auser, Basel (1992).

\bibitem[Wen]{Wen} H. Wente: An existence theorem for surfaces of
constant mean curvature, {\em J. Math. Anal. Appl.} {\bf 26} (1969),
318--344.

\bibitem[W]{W} T. J. Willmore: Total Curvature in Riemannian
Geometry, John Wiley \& Sons, New York (1982).

\bibitem[Z]{Z}M. Zhu: Harmonic maps from degenerating
Riemann surfaces, {\em Math. Z.} {\bf 264} (2010), 63--85.


\end{thebibliography}}
\end{document}